\documentclass[11pt]{article}
\usepackage{fullpage}
\usepackage{makeidx}

\newenvironment{@abssec}[1]{%
     \if@twocolumn
       \section*{#1}%
     \else
       \vspace{.05in}\footnotesize
       \parindent .2in
         {\bfseries #1. }\ignorespaces
     \fi}
     {\if@twocolumn\else\par\vspace{.1in}\fi}
\newenvironment{keywords}{\begin{@abssec}{Key words}}{\end{@abssec}}
\newenvironment{AMS}{\begin{@abssec}{AMS subject classification}}{\end{@abssec}}
\newtheorem{prop}{Proposition}
\newtheorem{thm}[prop]{Theorem}

\newtheorem{corol}[prop]{Corollary}

\def\diag{\mathop{\rm diag}\nolimits}

\def\eqbd{\mathop{{:}{=}}}

\def\openC{{\rm C\kern-.48em\vrule width.06em height.6em depth-.02em 
                 \kern.48em}}
\def\openR{{{\rm I}\kern-.16em {\rm R}}}
\def\openZ{{{\rm Z}\kern-.28em{\rm Z}}}
\def\sZZ{{{\scriptstyle\rm Z}\kern-.24em{\scriptstyle\rm Z}}}
\def\openT{{{\rm T}\kern-.42em {\rm T}}}
\def\openH{{{\rm I}\kern-.16em {\rm H}}}
\def\openK{{{\rm I}\kern-.16em {\rm K}}}
\def\openL{{{\rm I}\kern-.16em {\rm L}}}
\def\openM{{{\rm I}\kern-.16em {\rm M}}}
\def\openN{{{\rm I}\kern-.16em {\rm N}}}
\def\openP{{{\rm I}\kern-.16em {\rm P}}}

\def\eqbd{\mathop{{:}{=}}}

\def\sC{{\rm C\kern-.38em\vrule width.06em height.45em depth-.02em 
                 \kern.3em}}

\let\R\openR
\let\Z\openZ
\def\eop{\hfill
        {\ \vbox{\hrule\hbox{\vrule height1.3ex\hskip0.8ex\vrule}\hrule}}
        \vskip 0.3cm \par}

\def\belowrightarrow#1{{{{}\over\ #1\ }\kern-1.1em\to}}
\def\l2{{L_2}}

\def\Mu{{{\cal M}}}

\begin{document}

\title{The inverse eigenvalue problem for symmetric 
anti-bidiagonal matrices}
\author{ Olga Holtz \\ Department of Mathematics \\
University of California \\
Berkeley, California 94720 USA }
\date{\today}
\maketitle



\begin{abstract} The inverse eigenvalue problem for real symmetric
matrices of the form 
$$ \left [ \begin{array}{ccccccccc} 0 & 0 & 0 & & \cdots & & 0 & 0 & * \\
0 & 0 & 0 & & \cdots  & & 0 & * & * \\
0 & 0 & 0 & & \cdots  & & * & * & 0 \\
  &&  & &&  \cdot & & & \\
\vdots & \vdots & \vdots & & \cdot & & \vdots & \vdots & \vdots \\
 & & & \cdot &&   & & & \\
0 & 0 & * & & \cdots & & 0 & 0 & 0 \\
0 & * & * &  & \cdots & & 0 & 0 & 0 \\
{*} & {*} & 0  & & \cdots &  & 0 & 0 & 0 \end{array}  \right] $$
is solved. The solution is shown to be unique. The problem is also shown 
to be equivalent to the inverse eigenvalue problem for a certain subclass 
of Jacobi matrices. 
\end{abstract}

\section{Introduction}
Ths goal of this paper is to characterize completely the spectra
of real symmetric {\em anti-bidiagonal\/} matrices, i.e., matrices of the 
form  \begin{equation}
A= \left [ \begin{array}{ccccc} 0 & 0 & \cdots & 0 & a_n \\
0 & 0 & \cdots & a_{n-2} & a_{n-1} \\ \vdots & \vdots & \cdot & \vdots &
\vdots  \\0 & a_{n-2} & \cdots & 0 & 0 \\ a_{n} & a_{n-1} & \cdots & 0 & 0 
\end{array}  \right], \quad a_1, \ldots, a_n\in \R. 
\label{mainform}
\end{equation}

This work is motivated by the author's ongoing work on the nonnegative
inverse eigenvalue problem and is inspired by well-known results on
Jacobi matrices due to Hochstadt~\cite{Ho1},~\cite{Ho2}, 
Hald~\cite{Ha}, Gray and Wilson~\cite{GW}, as well as by the classical
 connection between the
Jacobi matrices and orthogonal polynomials (see, e.g.,~\cite[p.~267]{AAR}).
 
The blanket assumption for the rest of the paper is that all $a_j$ 
are positive. This restriction is clearly unimportant, since the sign of any 
$a_j$, $j>1$, can be changed using a unitary similarity of the form
$$ \diag (\varepsilon_1, \cdots, \varepsilon_n), \qquad \varepsilon_j=\pm 1, $$
and the problem for $a_1<0$ can be solved by switching from $A$ to $-A$.
The assumption $a_j>0$, $j=1,\ldots, n$, is however just right to guarantee
uniqueness of a matrix that realizes a given $n$-tuple as its spectrum.

\section{Definitions and notation}

Notation used in the paper is rather standard. The spectrum of a matrix
$A$ is denoted by $\sigma(A)$. A submatrix of $A$ with rows indexed by 
an increasing sequence $\alpha$ and columns indexed by another
sequence $\beta$ is denoted by $A(\alpha,\beta)$. For simplicity, 
a principal submatrix of $A$ with rows and columns indexed by $\alpha$ 
is denoted by $A(\alpha)$. (A typical choice 
for such an $\alpha$ will be $i{:}j$, the sequence of consecutive integers 
$i$ through $j$.) The size of a sequence $\alpha$ is denoted by $\# \alpha$.
If $\# \alpha=\# \beta$, then $\det A(\alpha,\beta)$ is denoted by
$A[\alpha,\beta]$; $\det A(\alpha)$ is denoted by $A[\alpha]$.
the elementary symmetric functions of an $n$-tuple $\Lambda$ are denoted
as $\sigma_j(\Lambda)$. Thus
$$  \sigma_1(\Lambda)\eqbd \sum_{j=1}^n \lambda_j, \qquad \sigma_2(\Lambda)\eqbd
\sum_{i<j} \lambda_i \lambda_j, \qquad etc.  $$
The term  anti-bidiagonal matrix was already introduced. Other
requisite definitions are listed next.

A {\em Jacobi matrix\/} is a 
tridiagonal matrix with positive codiagonal entries. 

A {\em sign-regular
matrix of class $d\leq n$ with signature sequence\/} $\varepsilon_1$, 
$\ldots$, $\varepsilon_d$, 
where $\varepsilon_j=\pm 1$ for all $j$, is a matrix satisfying
$$ \varepsilon_j A[\alpha] \geq 0 \quad {\rm whenever}\;\; \# \alpha=j, \quad
j=1, \ldots, d.  $$
If in addition all minors of order at most $d$ are nonzero, the
matrix is called {\em strictly sign-regular\/}. Finally, 
if a certain power of a sign-regular matrix of class $d$ is
strictly sign-regular, then the matrix is called sign-regular 
of class $d^+$. A particular case of strict sign regularity is {\em total
positivity\/} when all minors of a matrix are positive.
  
A sequence $\mu_1<\cdots <\mu_k$ is said to {\em interlace\/} a sequence
$\lambda_1<\ldots < \lambda_{k+1}$ if
$$  \lambda_1< \mu_1 < \lambda_2 <\mu_2 <\cdots < \mu_k < \lambda_k.  $$

\section{Results}

The following theorem is the main result of this paper.

\begin{thm} A real $n$-tuple $\Lambda$ can be realized as the spectrum 
of an anti-bidiagonal matrix~(\ref{mainform}) with all $a_j$ positive
if and only if $\Lambda=(\lambda_1, \ldots,\lambda_n)$ where
\begin{equation}
 \lambda_1>-\lambda_2>\lambda_3 >\cdots > (-1)^{n-1} \lambda_n>0. 
\label{spec} \end{equation}
The realizing matrix is necessarily unique. 
\end{thm}

\noindent {\bf Proof.\/\/} Necessity. Let $J$ denote the antidiagonal
unit matrix
$$ J\eqbd \left[ \begin{array}{ccccc}  0 & 0 & \cdots & 0 & 1 \\
 0 & 0 & \cdots & 1 & 0 \\ \vdots & \vdots & \cdot & \vdots & \vdots \\ 
0 & 1 & \cdots & 0 & 0 \\  1 & 0 & \cdots & 0 & 0 \end{array} \right].$$
Note that $J$ is sign-regular of class $n$ with the signature sequence 
\begin{equation} 1, -1, -1, 1, 1, \cdots, (-1)^{\lceil n-1/2 \rceil}. 
\label{signs} \end{equation}
Next, note that $B\eqbd JA$ is a nonnegative bidiagonal matrix, hence all its
minors are nonnegative. Now, by the Cauchy-Binet formula
$$  A[\alpha]=(JB)[\alpha]=\sum_{\# \beta=\#\alpha} J[\alpha,\beta] 
B[\beta,\alpha]. $$
Since the only nonzero minors of $J$ are principal, we conclude that the 
 matrix $A$ is sign-regular of class $n$ with the same signature sequence~(\ref{signs}). Since $A^2$ is a positive definite Jacobi matrix, a high enough
power of $A^2$ is totally positive, hence $A$ is sign-regular of type $n^+$.

By a theorem of Gantmacher and Krein~\cite[p.~301]{GK}, the eigenvalues of $A$ 
therefore can be arranged to form a sequence with alternating signs
and strictly decreasing absolute values whose first element is positive,
i.e., the spectrum $\sigma(A)$ satisfies~(\ref{spec}).

Sufficiency. First reduce the inverse problem for anti-bidiagonal 
matrices to the inverse problem for certain Jacobi matrices.
Consider a matrix of the form~(\ref{mainform}). To stress its dependence
on $n$ parameters $a_1$ through $a_n$, let us denote it by $A_n$. 
The argument will involve the collection of all matrices $A_n$, $n\in \Z$, 
determined by a single sequence $a_1$, $a_2$, $\ldots$.
Denote the characteristic polynomial of $A_n$
by $p_n$: $$ p_n(\lambda)\eqbd  \det (\lambda I - A_n).$$
Expanding it by its first row yields
\begin{eqnarray}
&& p_n(\lambda)=\lambda p_{n-1}(\lambda)-a_n^2p_{n-2}(\lambda), \quad n\geq 2
 \label{recur1}  \\
&& p_0(\lambda)=1, \qquad p_1(\lambda)=\lambda-a_1,   \label{recur2} 
\end{eqnarray} 
since the matrix $A_{n-1}$ is similar to its reflection about the 
antidiagonal. 

This three-term recurrence relation~(\ref{recur1}) with initial
conditions~(\ref{recur2}) is also satisfied~(see, 
e.g.,~\cite[p.~267]{AAR} or check directly) by the characteristic 
polynomials of the Jacobi matrices 
\begin{equation}
B_n\eqbd \left[ \begin{array}{cccccc} a_1 & a_2 & 0 & \cdots & 0 & 0 \\
 a_2 & 0 & a_3 & \cdots & 0 & 0 \\
 0  & a_3 & 0 & \cdots & 0 & 0 \\
\vdots & \vdots & \vdots & \ddots & \vdots & \vdots \\
0 & 0 & 0 & \cdots & 0 & a_n \\
0 & 0 & 0 & \cdots & a_n & 0  \end{array} \right]
\label{specjacobi}
\end{equation}
if each of them is expanded by its last row.
Thus the inverse eigenvalue problem for anti-bidiagonal matrices $A_n$ is 
equivalent  to the inverse eigenvalue problem for Jacobi matrices $B_n$.

Now comes the crucial step in the proof. Consider expanding the characteristic
polynomials of matrices $B_n$ in the opposite order, i.e., starting from the
first row. Precisely, let us denote by $q_{n-j+1}$ the characteristic polynomial
of the principal submatrix $B_n(j{:}n)$, with $q_n=p_n$. The corresponding 
recurrence relation is  
\begin{eqnarray}
&& q_n(\lambda)=(\lambda-a_1) q_{n-1}(\lambda)-a_2^2q_{n-2}(\lambda), 
\label{newrecur1} \\
&& q_{n-j}(\lambda)=\lambda q_{n-j-1}(\lambda)-a_{j+2}^2 q_{n-j-2}(\lambda),
\quad j=1,\ldots, n-2, \label{newrecur2}\\ 
&& q_0(\lambda)=1, \qquad q_1(\lambda)=\lambda.   \label{newrecur3} 
\end{eqnarray} 

Let $\Lambda$ be an $n$-tuple satisfying~(\ref{spec}). Define
the polynomial $q_n$ as
$$ q_n(\lambda)\eqbd \prod_{j=1}^n (\lambda-\lambda_j) $$
and show that one can define polynomials $q_{n-j}$ for all $j=1, \ldots, n$
so as to meet the requirements~(\ref{newrecur1})--(\ref{newrecur3}).
To this end, first define
\begin{equation} a_1\eqbd \sigma_1(\Lambda), \qquad  q_{n-1}(\lambda)=
{(-1)^n q_n(-\lambda)-q_n(\lambda) \over 2a_1}. \label{qn-1} 
\end{equation}
Note that $a_1>0$ due to the properties of $\Lambda$ and
that the (monic) polynomial $q_{n-1}$ is even or odd depending on whether
$n-1$ is even or odd. Also note that the coefficient of $\lambda^{n-3}$
in $q_{n-1}$ is equal to
$$ {\sigma_3(\Lambda)\over a_1}={\sigma_3(\Lambda)\over \sigma_1(\Lambda)}<0.$$
On the other hand, the coefficient of $\lambda^{n-2}$ in $q_n(\lambda)$ is
$\sigma_2(\Lambda)<0$. Therefore, it remains to show that the quantity
${\sigma_3(\Lambda) \over \sigma_1(\lambda)}-\sigma_2(\Lambda) $
is positive, so $a_2$ can be defined as its (positive) square root:
$$ a_2\eqbd\sqrt{{\sigma_3(\Lambda) \over \sigma_1(\lambda)}-\sigma_2(\Lambda)}. $$ 

Indeed, let us prove that 
\begin{equation}\sigma_3 >\sigma_1 \sigma_2 \label{sigmas} \end{equation}
 by induction. The base case is $n=3$, where $\lambda_1>-\lambda_2>\lambda_3>0$.
Then~(\ref{sigmas}) reduces to the inequality
\begin{equation}  \left( {1\over \lambda_1} +{1\over \lambda_2}
+{1\over \lambda_3}  \right)(\lambda_1+\lambda_2+\lambda_3)>1.
\label{lambdas}
\end{equation}
Differentiating the left-hand side of~(\ref{lambdas}), one can check that it is an 
increasing function of $\lambda_1$ for $\lambda_1\geq -\lambda_2$.
Since the left-hand side is exactly $1$ when $\lambda_1=-\lambda_2$, this
proves~(\ref{lambdas}) and therefore proves~(\ref{sigmas}). If $n>3$, 
also notice that inequality~(\ref{sigmas}) turns into equality for
$\lambda_1=-\lambda_2$, so it remains to argue that the difference
$\sigma_3-\sigma_1\sigma_2$ is an increasing function of $\lambda_1$
for $\lambda_1\geq -\lambda_2$.  But this is indeed the case, 
as can be seen by considering symmetric functions of
the set $\Lambda'\eqbd -\lambda_2, \ldots, -\lambda_n$. Since
$$ \sigma_1(\Lambda)=\lambda_1-\sigma_1(\Lambda'), \quad
\sigma_2(\lambda)=-\lambda_1\sigma_1(\Lambda')+\sigma_2(\Lambda'),
\quad \sigma_3(\Lambda)=\lambda \sigma_2(\Lambda')-\sigma_3(\Lambda'),$$
the inequality~(\ref{sigmas}) amounts to
$$ \lambda_1^2 \sigma_1(\Lambda')-\lambda_1 \sigma_1^2(\Lambda') +\sigma_1(\Lambda')
\sigma_2(\Lambda')-\sigma_3(\Lambda')>0,  $$
and the derivative of the last left-hand side is positive, since
$\lambda_1\geq \sigma_1(\Lambda')$. This completes the proof of~(\ref{sigmas}).
Thus, $a_2$ is well-defined.

With these definitions in place, define $q_{n-2}$ from~(\ref{newrecur1}), i.e.,
let $$   q_{n-2}(\lambda)\eqbd -{q_n(\lambda)-(\lambda-a_1)q_{n-1}(\lambda) 
\over a_2^2}. $$ Note that $q_{n-2}$ is a monic polynomial and is odd or even 
(precisely, it has the same parity as its leading term). 

Now show that the roots of $q_{n-1}$ interlace
those of $q_n$ and the roots of $q_{n-2}$ interlace those of $q_{n-1}$.
Note that the polynomials $p_n(\lambda)$ and $(-1)^{n-1}p_n(-\lambda)$
have the same sign on the intervals
$$ ( -|\lambda_1|, -|\lambda_2|), \;\;  
 ( -|\lambda_3|, -|\lambda_4|), \;\; \ldots, \;\; ( |\lambda_2|, |\lambda_1|). $$
Moreover, the sequence of these signs is alternating. 
The polynomial $q_{n-1}$ defined by~(\ref{qn-1}) therefore has exactly
$n-1$ real zeros, each of them between two consecutive zeros of $q_n$.
The implication for the root interlacing of $q_{n-2}$ and $q_{n-1}$ is immediate
and is a standard argument on orthogonal polynomials 
(cf.~\cite[Section~5.4]{AAR}). Due to the root interlacing of $q_{n-1}$
and $q_n$ and due to~(\ref{newrecur1}), the values of $q_{n-2}$ at
the zeros of $q_{n-1}$ form an alternating sequence. Therefore, the roots
of $q_{n-2}$ interlace those of $q_{n-1}$.

The rest of the argument is quite straightforward. With $q_{n-j}$ and $q_{n-j-1}$
defined, one defines $q_{n-j-2}$ from~(\ref{newrecur2}) making sure that
$a_{j+2}^2$ is indeed positive, for each $j=1,\ldots, n-2$.  The resulting 
monic polynomials will have alternating parities and interlacing roots. 
The quantity $a_{j+2}^2$ is to be set equal to the difference between the
second elementary symmetric function $\sigma_2$ of the roots of $p_{n-j-1}$ 
and  the second elementary symmetric function of the roots of $p_{n-j}$.
With a slight abuse of notation, this may be denoted by
$$ a_{j+2}^2=\sigma_2(p_{n-j-1})-\sigma_2(p_{n-j}).$$  
The roots of either polynomial are symmetric about $0$, therefore,
the corresponding second elementary symmetric function is simply 
$$ (-1)\cdot \hbox{\rm the sum of squares of all positive roots.}$$ By the
interlacing property, the sum of squares for $p_{n-j}$ exceeds that 
for $p_{n-j-1}$, hence  $\sigma_2(p_{n-j-1})-\sigma_2(p_{n-j})>0$
and hence $a_{j+2}$ is well defined.

The argument also shows the uniqueness of the realizing matrix~(\ref{specjacobi}),
therefore the uniqueness of the realizing matrix~(\ref{mainform}), provided,
of course, that $a_j$ are chosen to be positive.  \eop

The following corollary was established in the course of the above proof.

\begin{corol} A real $n$-tuple $\Lambda$ can be realized as the spectrum
of a Jacobi matrix~({\ref{specjacobi}}) if and only if
$\Lambda=(\lambda_1, \ldots,\lambda_n)$ where
$$ \lambda_1>-\lambda_2>\lambda_3 >\cdots > (-1)^{n-1} \lambda_n>0.$$ 
The realizing matrix is necessarily unique.

\end{corol}

Finally, another simple consequence of Theorem~1 is the following result.

\begin{corol} Let $\Mu$ be a real positive $n$-tuple. Then there exists
a Jacobi matrix that realizes $\Mu$ as its spectrum and has an anti-bidiagonal 
symmetric square root of the form~(\ref{mainform}) with all $a_j$ positive.
\end{corol}

\noindent {\bf Proof.\/\/} Let the elements of $\Mu$ be ordered 
$\mu_1>\mu_2>\cdots> \mu_n(>0)$. Define 
$$ \lambda_j\eqbd (-1)^{j-1}\sqrt{\mu_j}, \quad j=1, \ldots, n, \qquad
\Lambda\eqbd(\lambda_j : j=1, \ldots, n). $$
Then $$  \lambda_1>-\lambda_2>\lambda_3>\cdots>(-1)^{n-1} \lambda_n >0. $$
By Theorem~1, there exists a symmetric anti-bidiagonal matrix $A$ with
spectrum $\sigma(A)=\Lambda$. But then $B\eqbd A^2$ is a Jacobi matrix with 
spectrum $\Mu$. \eop

\section*{Acknowledgments} I am grateful to Plamen Koev for pointing me to
results on sign-regular matrices and to Gautam Bharali for finding a 
miscalculation in an earlier proof of Theorem~1.


\begin{thebibliography}{111}

\bibitem{AAR} G.E. Andrews, R. Askey, R. Roy, {\sc Special functions.\/} 
Cambridge University Press, Cambridge, 1999.

\bibitem{GK} F. R. Gantmacher and M. G. Krein, {\sc Oszillazionsmatrizen,
Oszillazionskerne und kleine Schwingungen mechanischer Systeme.\/}
Akademie-Verlag, Berlin, 1960.

\bibitem{GW} L. J. Gray and D. G. Wilson, {\em Construction of a Jacobi
matrix from spectral data, \/} Linear Algebra Appl., {\bf 14} (1976), 131--134.

\bibitem{Ha} O. Hald, {\em Inverse eigenvalue problems for Jacobi matrices,\/}
Linear Algebra Appl., {\bf 14} (1976), 63--85. 

\bibitem{Ho1} H. Hochstadt, {\em On some inverse problems in matrix theory,\/}
Arch. Math. {\bf 18} (1967), 201--207. 

\bibitem{Ho2} H. Hochstadt, {\em On the construction of a Jacobi matrix from spectral data,\/} Linear Algebra Appl., {\bf 8} (1974), 435--446.

\end{thebibliography}
\end{document}